\documentstyle[verbatim,righttag]{amsart}

\numberwithin{equation}{section}

\newcommand{\ms}{\medskip}

\newcommand{\noi}{\noindent}
\newcommand{\ra}{\rightarrow}
\newcommand{\bea}{\begin{eqnarray}}
\newcommand{\eea}{\end{eqnarray}}
\newcommand{\gr}{Groenewold}
\newcommand{\vh}{Van Hove}
\newcommand{\vn}{von Neumann}
\newcommand{\q}{\cal Q}

\newcommand{\h}{\cal H}
\newcommand{\f}{\cal F}
\newcommand{\oo}{\cal O}

\newtheorem{thm}{Theorem}

\newtheorem{prop}[thm]{Proposition}

\theoremstyle{definition}
\newtheorem{defn}{Definition}
\newsymbol\ltimes 226E
\newsymbol\restriction 1316

\begin{document}

\title{On a Full Quantization of the Torus}

\author{Mark J. Gotay}

\address{Department of Mathematics \\University of Hawai`i
\\2565 The Mall \\ Honolulu, HI 96822 USA}

\email{gotay@@math.hawaii.edu}

\date{July 31, 1995 (dg-ga/9507005)}

\thanks{Partially supported by NSF grant DMS-9222241}
\maketitle

\begin{abstract}

I exhibit a prequantization of the torus which is actually a
``full'' quantization in the sense that a certain complete
set of classical observables is irreducibly represented. Thus
in this instance there is no \gr\,-\vh\ obstruction to
quantization.

\end{abstract}

%%%%%%%%%%%%%%%%%%%%%%%%%%%%%%%%%%%%%%%%%%%%%%%%%%%%%%%%%%%%%%%%%%%

\section{Introduction}

The prequantization procedure produces a faithful
representation of the  entire Poisson algebra of a
quantizable symplectic manifold \cite{k}. In general, these
prequantization representations are flawed physically; for
instance, the prequantization of $\Bbb R^{2n}$ with its
standard symplectic structure is not unitarily equivalent to
the Schr\"odinger representation. One usually remedies this
by imposing an irreducibility requirement. But there is
seemingly a price to be paid for irreducibility: one can no
longer quantize {\em all\/} classical observables, but rather
only proper subalgebras thereof.

This ``obstruction'' to quantization has been known since the
1940s. In a series of papers, \gr\ and later \vh\ showed that
it is impossible to quantize the entire Poisson algebra of
polynomials on $\Bbb R^{2n}$ in such a way that the Heisenberg
\,h(2$n$)\, subalgebra of inhomogeneous linear polynomials is
irreducibly represented \cite{gr,vh1,vh2}. Some of the
maximal subalgebras of polynomials that {\em can} be
consistently quantized subject to this irreducibility
requirement are the inhomogeneous quadratic polynomials, and
polynomials which are at most affine in the momenta or the
configurations. See \cite{c1,f,j} for further discussions of
this example. Recently, a similar phenomenon was observed for
$S^2$
\cite{g-g-h}. In this case it was shown that the maximal
subalgebra of the Poisson algebra of spherical harmonics that
can be consistently quantized while irreducibly representing
the \,u(2)\, subalgebra of spherical harmonics of degree at
most one is just this \,u(2)\, subalgebra itself.

Based on these results as well as general quantization
theory, it would seem reasonable to conjecture that a
``no-go'' theorem must always hold, to wit:

\begin{quote}

{\em It is impossible to quantize the entire Poisson algebra
of any
\linebreak symplectic manifold subject to an irreducibility
requirement.}

\end{quote}

To make this conjecture precise, I introduce some
terminology. Let
$(M,\omega)$ be a connected symplectic manifold, and let
$\cal O$ be a Poisson subalgebra of
$C^{\infty}(M).$

\begin{defn} A {\sl prequantization} of $\oo$ is a linear map
$\q$ from
$\oo$ to an algebra of (essentially) self-adjoint
operators\footnote{For the most part technical difficulties
with unbounded operators will be ignored, as they are not
essential for what follows.} on a Hilbert space $\h$ such
that

\begin{enumerate}
\item[({\em i\/})] $\cal Q\big(\{f,g\}\big) = 2\pi i\big[\cal
Q(f),\cal Q(g)\big],$
\end{enumerate}

\noindent where $\{\,,\,\}$ denotes the Poisson bracket and
$[\;,\,]$ the commutator.\footnote{I use units in which
Planck's constant $h = 1$.} If
$\oo$ contains the constant function 1, then also
\begin{enumerate}
\item[({\em ii\/})] $\cal Q(1) = I.$
\end{enumerate}
\end{defn}

As the nomenclature suggests, it appears necessary in
practice to supplement these conditions, often by requiring
that a certain subset $\cal F$ of observables be represented
irreducibly. In favorable circumstances one can take $\cal F$
to consist of the components of a momentum map associated to
a transitive Lie symmetry group. In the nonhomogeneous case,
the corresponding notion is a ``complete set of
observables.'' This is a set
$\cal F$ with the property that $\{f,g\}=0$ for all
$f \in \f$ implies that $g = \mbox{const}.$ This means that
the Hamiltonian vector fields of elements of $\cal F$ span
the tangent spaces to $M$ almost everywhere. I shall assume
that the linear span of $\cal F$ is finite-dimensional, in
which case
$\cal F$ is said to be {\sl finite}.

Similarly, a set of operators $\cal A$ on $\h$ is {\sl
irreducible} provided the only (bounded) operator which
commutes with each $A \in \cal A$ is a multiple of the
identity. Since irreducibility is the quantum analogue of
completeness, I make

\begin{defn} Let $\cal F \subset \oo$ be a finite complete
set of observables. A prequantization $\q$ of $\oo$ is a {\sl
quantization} of the pair
$(\oo,\cal F)$ provided
\begin{enumerate}
\item[({\em iii\/})] the corresponding operators
$\big\{\q(f)\,|\,f \in \f\big\}$ form an irreducible set.
\end{enumerate}
\noi If one can take $\oo = C^{\infty}(M)$, the quantization
is {\sl full.}\footnote{\,Without some sort of finiteness
condition on $\cal F$, it is often possible to find full
quantizations; indeed, it may happen that a prequantization
representation is itself irreducible. For example, Chernoff
has constructed irreducible prequantization representations
when $M$ is compact
\cite{c2}.}
\end{defn}

Thus there do not exist full quantizations of either
$\big(\Bbb R^{2n},\mbox{h}(2n)\big)$ or
$\big(S^2,\mbox{u}(2)\big)$. According to the conjecture
above, it should be impossible to fully quantize {\em any}
symplectic manifold. But here I show that this is {\em
false\/}: there exists a full quantization of the torus
$T^2$. The existence of this full quantization, which can be
obtained as a particular Kostant-Souriau prequantization, is
surprising, especially given the \gr-\vh\ obstruction to a
full quantization of its covering
$\Bbb R^{2}$. The dichotomy stems from the fact that, unlike
on
$\Bbb R^2$, any prequantum bundle on $T^2$ is nontrivial.
This twisting forces the wave functions to be {\em
quasi-}periodic, and these boundary conditions effectively
override the obstruction in the case of interest. That the
prequantization I consider may indeed provide a full
quantization is signalled by the observation that it has, in
a sense, the Schr\"odinger representation of the
Heisenberg algebra built in. It turns out that this
representation of the Heisenberg algebra is well known in
solid state physics, where it goes under the name ``{\em
kq}-representation'' \cite{z}.

%%%%%%%%%%%%%%%%%%%%%%%%%%%%%%%%%%%%%%%%%%%%%%%%%%%%%%%%%%%%%%%%%%%

\section{Prequantization}

It is convenient to realize the torus $T^2$ as $\Bbb R^2/\Bbb
Z^2$ with symplectic form

\[\omega = Ndx \wedge dy,\]

\noi where $N$ is a nonzero integer. Then one can identify
$C^{\infty}(T^2)$ with the space of doubly periodic functions
$f$ on the plane:
\begin{eqnarray*} f(x+m,y+n) = f(x,y),& & m,n \in \Bbb Z.
\end{eqnarray*}

For the prequantization of the torus, I follow
\cite[\S2.3]{k}. (See also \cite{m} for a gauge-theoretic
treatment.) Let
$L_N$ be a prequantum bundle over $T^2$ with Chern class
$N.$\footnote{Up to equivalence, there is exactly one line
bundle per value of $N$, but each $L_N$ carries a
two-parameter family of inequivalent connections. Below I fix
a particular connection.} To explicitly construct
it,\footnote{The associated principal circle bundle over
$T^2$ may be realized as
$\Gamma_N\backslash H_N(2),$ where $H_N(2)$ is the
(polarized) Heisenberg group consisting of all matrices of
the form

\[ \left(\begin{array}{ccc} 1 & a & -c/N \\ 0 & 1 & b \\ 0 &
0 & 1
\end{array} \right) \]

\noi with $a,b,c \in \Bbb R$, and $\Gamma_N$ is the subgroup
of $H_N(2)$ consisting of those matrices for which $a$, $b$
and $c$ are integers
\cite[\S2]{f-g-g}.} introduce the sets

\[U_- = \big\{(x,y) \in (-\delta,1-\delta) \times
[0,1]\big\}\;\;\;\;\mbox{ and }\;\;\;\;  U_+ = \big\{(x,y)
\in (\delta,1+\delta) \times [0,1]\big\}\]

\noi for $\delta > 0$ small. Identifying $y=0$ with $y=1$
gives a pair of cylinders, and then identifying $x$ with
$x+1$ for $x
\in (-\delta,\delta)$ pastes the cylinders together into a
torus. On $U_{\pm}$ define the connection potentials

\[\theta_{\pm} = Nx\,dy.\]

\noi On the overlap $U_+ \cap U_-$, we have

\[\theta_- - \theta_+ = \left\{ \begin{array}{cl} 0, & x \in
(\delta,1-\delta)\\ -N\,dy, & x \in (-\delta,\delta)
\end{array} \right.\]

\noi whence the transition functions are

\[c(x,y) = \left\{ \begin{array}{cl}  1, & x \in
(\delta,1-\delta)\\ e^{-2\pi i Ny}, & x \in (-\delta,\delta).
\end{array} \right.\]

The space $\Gamma(L_N)$ of smooth sections of $L_N$ may thus
be identified with the space of smooth `quasi-periodic'
complex-valued functions
$\phi$ on the plane:
\begin{eqnarray}
\phi(x+m,y+n) = e^{2\pi iNmy}\phi(x,y),& & m,n \in \Bbb Z.
\label{qp}
\end{eqnarray}

\noi The corresponding prequantum Hilbert space $\h_N$
consists of those
$\phi(x,y)$ satisfying \eqref{qp} which are square integrable
over $[0,1)^2.$

The prequantization map $\q_N$ is
\begin{equation}
\q_N(f) =\frac{1}{2\pi iN}\bigg(\frac{\partial f}{\partial
x}\Big(\frac{\partial}{\partial y} - 2\pi iNx\Big)
-\frac{\partial f}{\partial y}\frac{\partial}{\partial
x}\bigg) + f.
\label{pm}
\end{equation}

\noi These operators are essentially self-adjoint on
$\Gamma(L_N).$

Now define, for each $N \neq 0$, the complete
sets\footnote{Although in principle I consider only
real-valued observables, here it is convenient to use complex
notation.}

\[\cal F_N =
\{e^{\pm 2\pi iNx},e^{\pm 2\pi iNy}\}.\]

\noi Then \eqref{pm} gives
\bea
\q_N(e^{\pm 2\pi iNx}) & = & e^{\pm 2\pi iNx}\bigg(1 \mp 2\pi
iNx
\pm \frac{\partial}{\partial y}\bigg) \label{op1}
\\
\q_N(e^{\pm 2\pi iNy}) & = & e^{\pm 2\pi iNy}\bigg(1 \mp
\frac{\partial}{\partial x}
\bigg).
\label{op2}
\eea

\noi Note that $\q_N(e^{-2\pi iN x}) = \q_N(e^{2\pi iNx})^*$
and
$\q_N(e^{-2\pi iN y}) = \q_N(e^{2\pi iNy})^*$.

\ms

%%%%%%%%%%%%%%%%%%%%%%%%%%%%%%%%%%%%%%%%%%%%%%%%%%%%%%%%%%%%%%%%%%%

\section{The Go Theorem}

The key observation is that the operators
\begin{equation}
\hat X = \frac{1}{2\pi i}\left(\frac{\partial}{\partial y} -
2\pi iNx \right),
\;\;\;\; \hat Y = -\frac{1}{2\pi i}\frac{\partial}{\partial
x}, \;\;\;\; \hat Z = -\frac{N}{2\pi i}\,,
\label{hat}
\end{equation}

\noi the first two of which which appear in \eqref{pm},
provide a representation of the Heisenberg algebra

\[\mbox{h}(2) = \big\{(X,Y,Z) \in \Bbb R^3\, |\, [X,Y] = Z,
\;[X,Z] = [Y,Z]= 0\big\}\]

\noi on $\h_N$. For $|N| = 1$ this   is equivalent to the
Schr\"odinger representation, as I now show.

For the moment fix $N=1$. ($N= -1$ will work as well.) The
analysis is simplified by applying the Weil-Brezin-Zak
transform $W$
\cite[\S1.10]{f},\cite{z} to the data in the previous
section. Let $\cal S(\Bbb R)$ be the Schwartz space. Define a
linear map $W: \cal S(\Bbb R) \ra \Gamma(L_1)$ by

\[(W\psi)(x,y) = \sum_{k \in \Bbb Z}\psi(x+k)e^{-2\pi iky}\]

\noi with inverse

\[(W^{-1}\phi)(x) = \int_0^1\phi(x,y)\,dy.\]

\noi That $W$ extends to a unitary map of $L^2(\Bbb R)$ onto
$\cal H_1$ is readily checked, cf. \cite[Thm. 1.109]{f}.
Under this transform,
$\hat X$ and
$\hat Y$ go over to the operators $-x$ and $-\frac{1}{2\pi
i}\frac{d}{dx}$ on
$L^2(\Bbb R)$, respectively. Thus the
$N = 1$ prequantization carries the Schr\"odinger
representation of the Heisenberg algebra.

This leads one to suspect that $\q_1$ might yield a full
quantization of the torus. Indeed this is the case:
\begin{thm}[Go Theorem] The prequantization $\q_1$ gives a
full quantization of
\linebreak
$\big (C^{\infty}(T^2),\cal F_1\big )$. \label{go}
\end{thm}
\begin{pf} Since $\q_1$ is a prequantization, conditions
({\em i}\/) and ({\em ii}\/) are automatically satisfied.
Thus it is only necessary to verify ({\em iii\/}), i.e.,

\[\q_1 \cal F_1 = \big\{\q_1(e^{\pm 2\pi ix}),\q_1(e^{\pm
2\pi iy})\big\}\]

\noi form an irreducible set. Applying $W$ to \eqref{op1} and
\eqref{op2}, one has, as operators on
$\cal S(\Bbb R)$,
\bea
\big(\q_1(e^{\pm 2\pi ix})\psi\big)(x) & = & e^{\pm 2\pi
ix}(1 \mp 2\pi ix)\psi(x)
\nonumber
\\
\big(\q_1(e^{\pm 2\pi iy})\psi\big)(x) & = & \bigg(1
\mp \frac{d}{dx}\bigg)\psi(x \pm 1).
\nonumber
\eea

Suppose that $T$ is a bounded linear operator on $L^2(\Bbb
R)$ which commutes with both
$A :=\q_1(e^{2\pi ix})$, $B :=\q_1(e^{2\pi iy})$ and their
adjoints. Then $T$ must commute with
$A^*A = 1 + 4\pi^2x^2$  and $B^*B= 1-d\,^2/dx^2.$ So
$[T,x^2]= 0$ and
$[T,d\,^2/dx^2]= 0,$ whence $T$ must commute with
$[d\,^2/dx^2,x^2] = 4x\,d/dx + 2.$ Thus $T$ commutes with the
generators of the symplectic algebra \,sp($2,\Bbb R$) in the
metaplectic representation $\mu$. Since
$\cal S(\Bbb R)$ contains a dense subspace of analytic
vectors for $\mu$ (viz. the Hermite functions, cf.
\cite[\S4.3]{f}), it follows that $[T,\mu(\frak M)] = 0$ for
all $\frak M$ in some neighborhood of the identity in the
metaplectic group
\,Mp$(2,\Bbb R)$ and hence, as this group is connected, for
all $\frak M$ in \,Mp$(2,\Bbb R)$.

Although the metaplectic representation is reducible, the
subrepresentations
$\mu_e$ and $\mu_o$ on each invariant summand of $L^2(\Bbb R)
= L^2_{e}(\Bbb R)
\oplus L^2_{o}(\Bbb R)$ of even and odd functions are
irreducible
\cite[\S4.4]{f}. Writing $T = P_eT + P_oT$, where $P_e$ and
$P_o$ are the even and odd projectors, one has
\begin{equation} [P_eT,\mu(\frak M)] = P_e[T,\mu(\frak M)] +
[P_e,\mu(\frak M)]T = 0
\label{com}
\end{equation}

\noi for any $\frak M \in \mbox{Mp(2,}\Bbb R)$. It then
follows from the irreducibility of the subrepresentation
$\mu_e$  that
$P_eT = k_eP_e + RP_o$  for some constant $k_e$ and some
operator $R:L^2_{o}(\Bbb R) \ra L^2_{e}(\Bbb R)$.
Substituting this expression into \eqref{com} yields
$[RP_o,\mu(\frak M)] = 0$, and Schur's Lemma then implies
that $R$ is either an isomorphism or is zero. But $R$ cannot
be an isomorphism as the representations
$\mu_e$ and
$\mu_o$ are inequivalent \cite[Thm. 4.56]{f}. Thus $P_eT =
k_eP_e$. Similarly
$P_oT = k_oP_o,$ whence $T = k_eP_e + k_oP_o.$

But now a short calculation shows that $T$ commutes with

\[A-A^* = 2i(\sin 2\pi x - 2\pi x\cos 2\pi x)\]

\noi only if $k_e = k_o$. Thus $T$ is a multiple of the
identity, and so
$\{A,A^*,B,B^*\}$ is an irreducible set, as was to be shown.
\end{pf}

This proof is not valid for $|N| > 1$. The problem is that
the map $W: L^2(\Bbb R)
\ra \h_N$, which in order to maintain \eqref{qp} now takes
the form

\[(W\psi)(x,y) = \sum_{k \in \Bbb Z}\psi(x+k)e^{-2\pi iNky},\]

\noi is no longer onto; indeed, from this formula it is
apparent that the image of
$W$ consists of those functions in $\h_N$ which have period
${1/N}$ in $y$. Denote the subspace of all such functions by
$\cal P_N$. Since the operators
$\hat X$ and
$\hat Y$ preserve periodicity, $\Gamma(L_N) \cap \cal P_N$ is
invariant and it follows that the representation
\eqref{hat} of
\,h(2)\, on $\h_N$ is no longer irreducible. While these
facts do not {\em a priori\/} preclude the existence of a
full quantization for $|N| > 1$, they do yield the
following.

\begin{prop} For $|N| > 1$, $\q_N$ does {\em not\/}
represent $\f_N$ irreduciblly.
\end{prop}

\begin{pf} It suffices to observe that the operators
\eqref{op1} and
\eqref{op2} commute with the orthogonal projector $\h_N \ra
\cal P_N.$
\end{pf}

Thus $\q_N$ provides a full quantization of
$\big(C^{\infty}(T^2),\f_N\big)$ iff
$|N| = 1.$

\ms

%%%%%%%%%%%%%%%%%%%%%%%%%%%%%%%%%%%%%%%%%%%%%%%%%%%%%%%%%%%%%%%%%%%

\section{Discussion}

The $|N| = 1$ quantization of the torus is curious in several
respects. First, it does not mesh well with what one
intuitively expects on the basis of geometric quantization
theory \cite{k}. For one thing, it is not necessary to
introduce a polarization in order to obtain an acceptable
quantization. (Since the prequantum wave functions are
quasi-periodic -- or, more crudely, since $\h_1
\approx L^2(\Bbb R)$ -- they are in effect ``already
polarized.'') And if one does polarize
$T^2,$ then one is guaranteed to be able to consistently
quantize only a substantially smaller set of observables --
namely, those whose Hamiltonian vector fields preserve the
polarization -- but not necessarily any larger subalgebra. On
the other hand, it is known that the torus admits a strict
deformation quantization \cite{r}, so perhaps it is not
entirely unexpected that it admits a full quantization as
well.

A second point concerns the role of the Heisenberg group
\,H(2)\, in this example, especially as compared to $\Bbb
R^2$. In both cases, it acts transitively (and factors
through a translation group.) On $\Bbb R^2$, there is a
momentum mapping for this action, and it is natural to insist
that quantization provide an irreducible representation of
\,h(2). But on $T^2$, there is no momentum map; the torus is
``classically anomalous'' in this regard \cite{a}. So it does
not make sense to require that quantization produce a
representation of \,h(2)\, in this case. Nonetheless, there
is the representation
\eqref{hat} of the Heisenberg algebra on
$\h_N$ which, for $|N| = 1$, is irreducible. The reason one
can obtain a full quantization while irreducibly representing
the Heisenberg algebra in this instance is because $\hat X$
and
$\hat Y$ are not the quantizations of any observables; in
other words, the representation
\eqref{hat} does not arise via quantization. Thus in this
sense the quantization is also anomalous. On
$\Bbb R^2$ this is not the case: demanding that
$\q$ be
\,h(2)-equivariant is too stringent a requirement to allow
the quantization of every observable in the Schr\"odinger
representation.

As an aside, it is interesting to observe that if $x$ and $y$
{\em were} globally defined observables, then from \eqref{pm}
one would have
\begin{equation}
\q (x) = \frac{1}{2\pi iN}\frac{\partial}{\partial
y}\;\;\;\;\mbox{   and   }
\;\;\;\;\q (y) = -\frac{1}{2\pi
iN}\left(\frac{\partial}{\partial x} - 2\pi iNy\right).
\label{qs}
\end{equation}

\noi Of course the operators
$\partial/\partial y$ and $\partial/\partial x - 2\pi iNy$
are not defined on
$\Gamma(L_N)$,\footnote{In \cite{a-c-g} such operators are
called ``bad.''} although they do make sense on
$\Gamma(L_{-N})$ with the connection given by the potentials
$\bar\theta_{\pm} = Ny\,dx$. In fact, one has

\[\frac{\partial}{\partial y} =
\overline\nabla_{\frac{\partial}{\partial
y}}\;\;\;\;\mbox{and}\;\;\;\;\frac{\partial}{\partial x} -
2\pi iNy =
\overline\nabla_{\frac{\partial}{\partial x}}\]

\noi where $\overline\nabla$ is the aforesaid connection on
$L_{-N}$. This explains the remark of Asorey \cite[\S4]{a},
to the effect that `the quantum generators of translation
symmetries are given by multiples of
$-i\overline\nabla_j$ for $j = 1,2$.' Furthermore, if one
were on $\Bbb R^2$ rather than $T^2$, the operators
\eqref{qs} would be well defined and would commute with
$\hat X$ and
$\hat Y$, indicating that the representation \eqref{hat} of
\,h(2)\, is reducible. So it is apparent how the
nontriviality of the prequantum bundles removes this
obstruction to the irreducibility of this representation
of the  Heisenberg algebra.

It is perhaps also surprising that the irreducibility of the
representation
\eqref{hat} of \,h(2)\, on $\h_1$ apparently does not, in and
by itself, imply that
$\q_1$ is a full quantization of
$\big(C^{\infty}(T^2),\f_1\big)$. Indeed, that the prequantum
Hilbert space carries the Schr\"odinger representation only
seems to be a portent; the proof of the Go Theorem devolves
instead upon the properties of the metaplectic
representation. I do not know if this theorem can be proved
directly using the irreducibility of \,h(2). More generally,
I wonder to what extent the fullness of the prequantization
$\q_N$, for a given complete set $\f$, is correlated with the
irreducibility of the \,h(2)\, representation? For instance,
is it possible to obtain a full quantization using $\q_N$ for
$|N| > 1$? Or perhaps one could prove a no-go result in this
context, thereby strengthening Proposition 2?

Thirdly, the requirement that quantization irreducibly
represent a complete set
$\f$ typically leads to ``\vn\ rules'' for elements of $\f$
\cite{g-g-h}. Roughly speaking, these rules govern the extent
to which $\q$ preserves the multiplicative structure of the
Poisson algebra. In particular, they determine how
$\q\big(f^2\big)$ is related to $\q(f)^2$ for $f \in \f.$ For
$\Bbb R^{2n}$ and
$S^2$ this relation is relatively `tight.' But for $T^2$ both
$\q_N\big(f^2\big)$ and $\q_N(f)^2$ are completely determined
for any observable $f$ and any $N$ by the simple fact that
$\q_N$ is a prequantization; irreducibility is irrelevant.
Moreover, one sees from
\eqref{pm} that $\q_N\big(f^2\big)$ is a first order
differential operator whereas $\q_N(f)^2$ is of second order.
This indicates that, in general, one cannot expect
quantization to respect the classical multiplicative
structure to any significant degree.

Another difference between the quantization of $T^2$ and the
quantizations of
$\Bbb R^{2n}$ and
$S^2$ is that in the latter cases the Poisson subalgebras
generated by the complete sets $\cal F$ are finite
dimensional. But for the torus, even though
$\cal F_1$ is finite, the Poisson subalgebra it generates --
the trigonometric polynomials -- is not, and is actually
dense in
$C^{\infty}(T^2).$

An alternate approach to the quantization of the torus, which
touches upon some of the issues discussed here, and which
contains applications to the quantum Hall effect, is given in
the recent paper by Aldaya et al.
\cite{a-c-g}.

%%%%%%%%%%%%%%%%%%%%%%%%%%%%%%%%%%%%%%%%%%%%%%%%%%%%%%%%%%%%%%%%%%%

\section*{Acknowledgments}

I thank Victor Aldaya and Julio Guerrero for helpful
discussions on no-go theorems.
\bigskip

%%%%%%%%%%%%%%%%%%%%%%%%%%%%%%%%%%%%%%%%%%%%%%%%%%%%%%%%%%%%%%%%%%%

%%%%%%%%%%%%%%%%%%%%%%%%%%%%%%%%%%%%%%%%%%%%%%%%%%%%%%%%%%%%%%%%%%%

\end{document}